\def\sqr#1#2{{\vcenter{\vbox{\hrule  height.#2pt
	\hbox{\vrule width.#2pt height#1pt \kern#1pt \vrule width.#2pt}
	\hrule height.#2pt}}}}
\def\sq{\sqr55}        
\documentclass[11pt]{article}
\usepackage{amssymb,amsmath,amscd}
\newcommand{\stdspace}{\hskip 0.75em\ignorespaces}

\renewenvironment{abstract}{\par{\bf Abstract}\par}{\par}

\newcommand{\co}{\colon\thinspace}    
\newcommand{\nl}{\hfil\break}         
\newcommand{\ppar}{\par\vskip 8pt plus4pt minus4pt}

\let\fonote\footnote
\def\footnote#1{\fonote{\small #1}}

\def\thanks#1{\relax}
\def\theoremstyle#1{\relax}
\def\dedicatory#1{\relax}
\def\tableofcontents{\relax}
\def\maketitle{\relax}
\def\date#1{\relax}
\def\bibliographystyle#1{\relax}
\def\title#1{\def\thetitle{#1}}

\def\author#1{\def\theauthors{#1}}
\def\address#1{\def\theaddress{#1}}

\def\email#1{\def\theemail{#1}}
\def\url#1{\def\theurl{#1}}
\long\def\abstract#1\endabstract{\long\def\theabstract{#1}}

\def\secondaryclass#1{\def\thesecondaryclass{#1}}
\def\keywords#1{\def\thekeywords{#1}}

\def\subjclass#1{\def\theprimaryclass{#1}}
\newtheorem{conjecture}[equation]{Conjecture}
\newenvironment{Correspondence}{\ppar\bf Main Theorem\sl\ \ }{\ppar}

\theoremstyle{definition}

\theoremstyle{remark}
\newtheorem{remark}[equation]{Remark}

\numberwithin{equation}{section}

\newcommand{\abs}[1]{\lvert#1\rvert}

\newcommand{\T}{{\widetilde{T\mathbb{P}}{}^1}}

\newcommand{\deq}{\overset{\text{def}}{=}}

\begin{document}
%
%
\input gtoutput
\volumenumber{1}\papernumber{2}\volumeyear{1997}
\published{14 July 1997}\pagenumbers{9}{20}

\title{Ward's Solitons}
\shorttitle{Ward's solitons}

\author{Christopher Anand}  

\address{DPMMS, 16 Mill Lane\\Cambridge, CB2 1SB, UK}

\email{C.Anand@pmms.cam.ac.uk}
\let\theurl\relax

\abstract
Using the `Riemann Problem with zeros' method, Ward has
constructed exact solutions to a $(2+1)$--dimensional integrable Chiral
Model, which exhibit solitons with nontrivial scattering.
We give a correspondence between what we conjecture to be all pure soliton
solutions and
certain holomorphic vector bundles on a compact surface.
\endabstract


\secondaryclass{58F07}


\proposed{Simon Donaldson}\seconded{Frances Kirwan, Peter Kronheimer}
\received{4 December 1996}\revised{16 May 1997}

\maketitlepage
%
\section{Introduction}

Nonlinear equations admitting soliton solutions in $3$--dimensional
space-time have been studied recently both numerically and
analytically.  See \cite{Su} and \cite{Wa2} for a discussion of
solitons in planar models.

In this paper, we  study an integrable model
introduced by Ward which is remarkable in that it possesses
interacting soliton solutions of finite energy \cite{{Su}, {Wa2}, {Io}}.
This \( \mathrm{SU}(N) \) chiral model with torsion term
may be obtained by dimensional reduction and gauge fixing from the \(
(2+2) \) Yang--Mills equations \cite{{Wa2}} or more directly from the
\( (2+1) \) Bogomolny equations.  Static solutions of the model
correspond to harmonic maps of \( \mathbb{R}^{2}\to{\mathrm{U}(N)} \)
which extend analytically to \( \mathbb{S}^{2} \) iff they have
finite energy.

The basic equations of Ward are
\begin{equation}\label{A-1}
	\frac{\partial}{\partial t}\bigl(J^{-1}\frac{\partial}{\partial t}J
	\bigr) - \frac{\partial}{\partial x}\bigr(J^{-1}
	\frac{\partial}{\partial x}J\bigr) -\frac{\partial}{\partial
	y}\bigl(J^{-1}\frac{\partial}{\partial y}J\bigr)+\bigl[J^{-1}
	\frac{\partial}{\partial y}J,J^{-1}\frac{\partial}{\partial
	t}J\bigr] =0
\end{equation}
where \( J\co \mathbb{R}^{3}\to\mathrm{SU}(N) \).  To this equation Ward
 added the boundary condition:
\begin{equation}\label{A-2}
	J(r,\theta,t)=\mathbb{I}+\frac{1}{r}
 J_{1}(\theta)+\mathcal{O}\left(\frac{1}{r^{2}}\right)
 	\quad\text{ as }\quad r\to\infty;
\end{equation}
we will assume \( J_{1}(\theta) \) is continuous.
Ward showed that analytic solutions to \eqref{A-1} correspond to
doubly-framed holomorphic bundles on the open surface \( T\mathbb{P}^{1} \).
We will show that a neccessary and sufficient condition for the
bundle to extend to the compactification \( \T \), the second
Hirzebruch surface is that \( J \) be analytic and that
the operator
 \begin{equation}
	  \frac{d}{du}+\tfrac{1}{2}
	 (1+\cos\theta)\iota^{*}\left(J^{-1}\frac{\partial }{\partial
	 x}J\right)
	 +\tfrac{1}{2}\sin\theta\,
\iota^{*}\left(J^{-1}\left(\frac{\partial }{\partial y}
	   +\frac{\partial}{\partial t}\right) J\right)
 	\label{A-3}
 \end{equation}
have null monodromy around  \( u\in\mathbb{R}\cup\{\infty\} \), where
  \begin{equation}
  	\iota(u)\deq (\cos\theta\, u+x_{0},\sin\theta\,u+y_{0},0),
  	\label{uline}
  \end{equation}
for all \( x_{0},y_{0}\in\mathbb{R} \) and \( \theta\in\mathbb{S}^{1} \),
i.e.\ for all lines in \( \mathbb{R}^{2} \).
    There is
some evidence that our techniques can be applied to the case of
nonanalytic solutions, but we will not do so here.  We also leave
open the question as to whether these are all the pure soliton
solutions.

Before going on, consider the null monodromy of \eqref{A-3} in the \(
\mathrm{U}(1) \) case, i.e.  for the usual d'Alembert equation.  Let
\( j=\log J \) be some logarithm of a solution.  The monodromy of
\eqref{A-3} becomes
\begin{displaymath}
	\int_{-\infty}^{\infty} \left[ (1+\cos\theta)\iota^{*}j_{x}
	+\sin\theta(\iota^{*}j_{y}+\iota^{*}j_{t})\right] du=0
\end{displaymath}
where \( j_{x}= \frac{\partial j}{\partial x}\), etc.
The fundamental theorem of calculus and the boundary condition
\eqref{A-2} imply
\begin{displaymath}
	\int_{-\infty}^{\infty}
\cos\theta\,\iota^{*}j_{x}+\sin\theta\,\iota^{*}j_{y}\
	du=0.
\end{displaymath}
Combining the two integrals with \( \theta=\theta_0 \) and
\( \theta = \theta_0+\pi \), we obtain
\begin{displaymath}
	0=\int_{-\infty}^{\infty}\sin\theta_0\,\iota^{*}j_{t}\ du=\sin\theta_0\
	  \frac{\partial }{\partial t}\int_{-\infty}^{\infty} \iota^{*} j\ du
\end{displaymath}
and
\begin{displaymath}
	0=\int_{-\infty}^{\infty}\iota^{*}j_{x}\ du.
\end{displaymath}
The first statement is that the Radon transform of \( j \) on a
space-plane is independent of time, and hence \( j \) is a harmonic
function.  Since \( j \) is also bounded (a result of \eqref{A-2}) it must be
constant.  This provides some support for the idea that \eqref{A-3}
has null monodromy for pure soliton solutions only.

We explain (in \S4) how
the boundary conditions can be interpreted in terms of the extension of
the holomorphic bundle to the fibrewise compactification (\( \T \))
when \( J \) satisfies \eqref{A-2} and \eqref{A-3} has null monodromy,
and to infinite points for fibres not above
the equator in \( \mathbb{P}^{1} \) (i.e.\ \(
\{\lambda\in\mathbb{C}\cup\{\infty\}:\abs{\lambda}\ne 1\} \)), when
\( J \) satisfies   \eqref{A-2} alone.

When \eqref{A-3} does have null monodromy, Serre's GAGA principle
tells us that the associated bundles are algebraic.  This explains the
algebraic nature of the solutions constructed so far, and was a strong
motivation for proving the main theorem.

\begin{Correspondence}\label{correspondence}  \sl
 There are bijections between the sets of 

\ppar
\leftskip 25pt
\noindent\llap{{\bf1\rm)}~~}%
analytic solutions \( J \) of \eqref{A-1} satisfying \eqref{A-2} for
which \eqref{A-3} has null monodromy; and 

\noindent\llap{{\bf2\rm)}~~}%
holomorphic rank \( N \) bundles \( \mathcal{V}\to T\mathbb{P}^{1}
 \) which are real in the sense that they admit a lift
 \begin{equation}\label{real structure}
 \begin{matrix}
 \begin{CD}
	 		\mathcal{V}@>\tilde{\sigma}>>\mathcal{V}\\
	 		@VVV@VVV\\
	 		T\mathbb{P}^{1}@>\sigma>>T\mathbb{P}^{1}
	 	\end{CD}  \quad
 	 	\begin{matrix}
	 		\text{of the}\\
	 		\text{ antiholomorphic}\\
	 		\text{ involution}
	 	\end{matrix}   \qquad
   	\begin{aligned}
	 		\sigma^{*}\lambda&=1/\bar{\lambda}\\
	 		\sigma^{*}\eta&=-\bar{\lambda}^{-2}\bar{\eta}
	 	\end{aligned}
 	\end{matrix}
 \end{equation}
 (where \( \lambda \) and \( \eta \) are standard base and fibre
 coordinates of \( T\mathbb{C}\subset T\mathbb{P}^{1} \))
 and which extend
 to bundles on the singular quadric cone \( T\mathbb{P}^{1}\cup
 \{\infty\} \), such that restricted to real sections (sections
 invariant under the real structure) \( \mathcal{V} \) is trivial, and
 restricted to the compactified tangent
 planes \(
 T_{\lambda}\mathbb{P}^{1}\cup \{\infty\} \) for \( \abs{\lambda}=1 \),
 \( \mathcal{V} \) is trivial,
 with a fixed, real framing.\par

 \end{Correspondence}

	\begin{remark} \rm The null monodromy of \eqref{A-3} makes sense
	for initial conditions on a space-plane \( \{t=t_{0}\} \).  It
	follows from the proof that the initial value problem with
	null-monodromy initial conditions has an analytic solution
	extending forward and backward to all time, i.e.  it cannot blow
	up in finite time.
\end{remark}

\subsection*{Construction of solutions}
There are currently three methods of solving this system.  The first
method of Ward was to give a twistor correspondence between solutions
of \eqref{A-1} and holomorphic bundles on \( T\mathbb{P}^{1} \), the
holomorphic tangent space to the complex projective line.  This led to
the construction of noninteracting soliton solutions.  Thereafter,
numerical simulations of these solutions by Sutcliffe led to his
discovery of interacting soliton solutions.  Exact solutions with two
interacting solitons were then constructed by Ward using a
Zakharov--Shabat procedure.  Using this procedure, more general
solutions were constructed by Ioannidou concurrently with the present
work.  In a future paper, we will present a closed-form expression for
all solutions satisfying \eqref{A-1}, \eqref{A-2} with null
\eqref{A-3} monodromy, including all known exact soliton solutions.
This will build on the monad-theoretic work in \cite{An2}.

\subsection*{Acknowledgements}
I am grateful to Sir Michael Atiyah, Riyushi Goto,  Partha Guha,
Nigel Hitchin, Mario Micallef, Paul Norbury, Richard Palais,
John Rawnsley, the referee and especially
Richard Ward for helpful discussions and advice.
This research was supported by an NSERC postdoctoral fellowship.

\section{Zero Curvature and the Bogomolny equations}\label{reduction}
Ward's equations are not a reduction in the sense
of dimensional reduction.  We obtain them from the Bogomolny
equations by fixing a gauge.

On $\mathbb{R}^{2+1}$, the Bogomolny equations for a connection
$\nabla=d+A$ and a Higgs
field (section of the adjoint bundle) $\Phi$ are
\begin{subequations}\label{Bogeqs}
		\begin{align}
		   -\nabla_{t}\Phi &= [\nabla_{x},\nabla_{y}]\\
			\nabla_{x}\Phi &=
[\nabla_{y},\nabla_{t}]\label{Bogeqsb}\\
			\nabla_{y}\Phi &= [\nabla_{t},\nabla_{x}].
		\end{align}
\end{subequations}

They are completely integrable, and can be written in the form
\begin{equation}
	[\nabla_{\bar z}+\frac{i\lambda}{2}\nabla_{t}-\frac{\lambda}{2}\Phi,
	\nabla_{z}-\frac{i}{2\lambda}\nabla_{t}-\frac{1}{2\lambda}\Phi] =0
	\quad\text{ for all } \lambda\in \mathbb{C}^{*}.
	\label{zccond}
\end{equation}
When \( \abs{\lambda}=1 \) this is the curvature for an underlying
connection on a family of planes.  Integrating it, we obtain a circle
of special gauges in which
\begin{align*}
	\Phi&=\Re\lambda\, A_{x}+\Im\lambda\, A_{y}\\
	A_{t}&=\Im\lambda\, A_{x}-\Re\lambda\, A_{y}.
\end{align*}
Ward's equations are equations for the gauge
transformation from the $\lambda=-1$ gauge to the $\lambda=1$ gauge.
We will call the $\lambda=-1$ gauge the {\it standard gauge}.

If $J$ is the gauge transformation, \eqref{Bogeqsb}
is Ward's equation \eqref{A-1}, and in the standard gauge,
\( \nabla=d+A \) and \( \Phi \) are
\begin{equation}
	\begin{aligned}
		-A_{x}=\Phi  &=\frac{1}{2} J^{-1}\partial_{x}J \\
		A_{y}=A_{t}  &=\frac{1}{2} J^{-1}
		    \bigl(\partial_{y}+\partial_{t}\bigr) J.
	\end{aligned}
	\label{in-1gauge}
\end{equation}
Conversely, given \( J \), we can form \( (\nabla,\Phi) \) in this
way.  Moreover, if \( J \) satisfies
\eqref{A-2} and has null \eqref{A-3} monodromy, the resulting map \(
J(z,t=0,\lambda)\co \mathbb{R}^{2}\times\mathbb{S}^{1}\to \mathrm{SU}(N) \)
extends to a based
map \( \mathbb{S}^{2}\times\mathbb{S}^{1}\to \mathrm{SU}(N) \).  This
associates a topological charge in \(
\pi^{3}(\mathrm{SU}(N))=\mathbb{Z} \) to any such solution \( J \).

\begin{conjecture} \sl
	This topological degree can be defined for all finite-energy
solutions,
	and is equal to the energy minus the effect of Lorentz boosting,
	internal spinning and radiation.
\end{conjecture}

\section{Twistor constructions of Ward and Hitchin}
Hitchin showed that the set of oriented geodesics on an
odd-dimensional real manifold has a complex structure (\cite{Hi}).  In
particular,
the set of
lines in $\mathbb{R}^{3}$ is isomorphic as a complex manifold to the
holomorphic tangent bundle of the complex projective line.
Using this equivalence he shows that solutions to the Bogomolny equations
correspond to holomorphic bundles on $T\mathbb{P}^{1}$.

Very briefly, given a solution $(\nabla,\Phi)$ to the Bogomolny
equations, one associates to a line the vector space of covariant
constant frames of the modified connection $\nabla -i\Phi$ on
 the line.  This is a complex bundle.  The operator $\nabla
_{\bar\eta}$ where $\eta$ represents a holomorphic fibre coordinate on
$T\mathbb{P}^{1}$ commutes with $\nabla -i\Phi$, and hence
descends to a $\bar{\partial}$--operator on the bundle.

The key point is the commuting of the two operators and after a
recombination, this can be written as a zero curvature condition.
See \cite{Hi} for a full account.

\section{The holomorphic bundle}
Given a solution $J$, let
$(\nabla=d+A,\Phi)$ be the solution to the Bogomolny equations,
in the standard gauge, as in \eqref{in-1gauge}.
The extension to the
compactification requires one argument near the equator
(\(\abs\lambda=1\)) (which
requires null \eqref{A-3} monodromy)
and another on the open hemispheres.

\setcounter{subsection}{\value{equation}}
\addtocounter{equation}{1}
\subsection{Away from the equator}
Consider
the $z$--plane, $\{t=0\}$, and the `projection':
 \begin{equation}
	 T\mathbb{S}^{2}\hookrightarrow\mathbb{R}^3 \times \mathbb{S}^2
	\to\mathbb{R}^2 \times\mathbb{S}^2
 \end{equation}
onto this plane.

The zero curvature connection has a
characteristic direction in this plane, and the appropriate linear
combination of the operators in \eqref{zccond} gives the
$\bar{\partial}$--operator for a rank $N$
bundle $\mathcal{V}_{\lambda}\to T_{\lambda}\mathbb{P}^{1}$:
\begin{equation}\label{op}
 	\overline{\nabla}_{\eta} \deq
  	(1+\lambda^{2})\partial_{x}+i(1-\lambda^{2})\partial_{y}
 	+(1+\lambda)^{2}A_{x} +i(1-\lambda^{2})(A_{y}+A_{t}).
\end{equation}
The kernel of this operator is the set of
 holomorphic sections of a bundle with
respect to the complex variable
\begin{equation}\label{comvar}
	\eta=\frac{z-\lambda^{2}\bar{z}}{1-\lambda\bar{\lambda}}.
\end{equation}
Together with $\partial_{\bar{\lambda}}$, this defines an operator
$$\overline{\nabla}\co \operatorname{gl}(\mathbb{C}^{N})\to
\operatorname{gl}\biggl(\mathbb{C}^{N}
\otimes T^{(0,1)}\bigl\{\abs{\lambda}<1,\eta \in
\mathbb{C}\bigr\}\biggr).$$
  Since $\overline{\nabla}_{\eta}$ depends holomorphically
on $\lambda$, $\overline{\nabla}^{2}=0$.  Under the assumption that
$J\in C^{1}(\mathbb{R}^{3})$ plus boundary conditions \eqref{A-2},
$\overline{\nabla}$ will be continuous on
$\{{\abs{\lambda}<1,\eta\in\mathbb{C}}\}$ which we identify with
$\{{\abs{\lambda}<1,z\in\mathbb{R}^{2}}\}$.
Near $z=\infty$
$$-\overline{\eta}^{2}\overline{\nabla}_{\eta}=
\partial_{1/\overline{\eta}}+C_{1}(\lambda,\theta) r^{2}A_{x}+
C_{2}(\lambda,\theta)r^{2}(A_{y}+A_{t})$$
where $z=re^{i\theta}$, and the functions $C_{1}$ and $C_{2}$ are
bounded in $\theta$ for each fixed $\lambda$, i.e.\ they are
polynomials in $\sin\theta$ and $\cos\theta$.  The boundary
conditions \eqref{A-2} for $J$ imply
\begin{equation}
	\begin{split}
		A_{x}&=J^{-1}(\cos\theta\partial_{r}+\frac{\sin\theta}{2ir}
		\partial_{\theta})J\\
		&=1/r^{2} A^{\prime}_{x}(1/r,\theta,t)
		\label{atzinfty}
	\end{split}
\end{equation}
where $A'_{x}$ is continuous near $z=\infty$, and similarly for
$A_{y}$ and $A_{t}$.  Hence $\overline{\nabla}_{\eta}$ is
continuous with a bounded singularity at $z=\infty$.

This implies that the coefficient is \(
L^{p}_{\text{loc}}(\mathbb{S}^{2}) \) for
\( 0<p\le\infty \) which is sufficient to show that iterating
convolution with the Cauchy kernel produces local holomorphic
gauges.  Since the data vary holomorphically in \( \lambda \), the
gauges can be used to define a holomorphic structure on \(
\mathcal{V}\to (\T \cap \{\abs{\lambda}<1\}) \).

\begin{remark} \rm
	The extension to the compactified nonequatorial fibres does not
	require the null \eqref{A-3} monodromy, and thus gives a
	necessary but not necessarily sufficient condition for a bundle to
	represent a solution satisfying the weak boundary condition.
\end{remark}

\setcounter{subsection}{\value{equation}}
\addtocounter{equation}{1}
\subsection{Null monodromy and the equator}
In the last section, we found a `\( \bar\partial \)--operator'
hidden in the zero curvature condition \eqref{zccond}.  Away from the
poles, we can make a different recombination of the operators, which
on the equator can be written in the manifestly real form
 \begin{equation}
	 \cos \theta\, \frac{\partial}{\partial
x}+\sin\theta\,\frac{\partial}{\partial y}
	 +\tfrac{1}{2}(1+\cos\theta)J^{-1}\frac{\partial }{\partial x}J
	 +\tfrac{1}{2}\sin\theta\ J^{-1}\left(\frac{\partial }{\partial y}
	   +\frac{\partial}{\partial t}\right) J.
 	\label{Du-Phi}
 \end{equation}
Under the assumption that \( J \) is analytic, this represents an
\( \mathbb{S}^{1}\times\mathbb{R} \)--family of first order ODEs on the
line which vary analytically with the parameter
\(\theta\in\mathbb{R} \).  The boundary
condition \eqref{A-2} implies that the functions
\( r^{2}J^{-1}\frac{\partial }{\partial x}J \) and
\( r^{2}J^{-1}\frac{\partial }{\partial y}J \) are bounded on
\( \mathbb{R}^{2} \), which means that
\( J^{-1}dJ \) has at worst a bounded discontinuity on \(
\mathbb{S}^{2} \), the conformal compactification of a space plane.
Since  the \( L^{1} \) norm is the natural norm in this context,
 we can convert all the
integrals on infinite lines to integrals over compact circles through
\( \infty\in\mathbb{S}^{2} \).  It follows that the coefficients vary
continuously in  \( L^{1} \) with the choice of line, and it makes
sense, given \( \theta \), to solve the whole family of ODEs on
parallel lines giving a function \( \mathbb{S}^{2}\to\mathrm{U}(N) \),
which is continuous at \( \infty \) iff \eqref{A-3} has null monodromy.

The result is an analytic map from \( \{\theta\in\mathbb{S}^{1}\} \) to
 \( C^{0}(\mathbb{S}^{2},\mathrm{U}(N)) \).
By analytic, we mean that it can be expanded in local
power series in \( \theta \) with coefficients in
\( C^{0}(\mathbb{S}^{2},\mathrm{U}(N)) \), which converge in some
neighbourhood with respect to the \( L^{\infty} \) norm
 (measured pointwise by geodesic distance from the unit in \( \mathrm{U}(N)
\)).
 This follows
from the fact that the operator \eqref{Du-Phi} is analytic in
\( \theta \) and hence has a power series which (in particular)
converges in the
\( L^{1} \) norm, and the integration map which solves the initial
value problem is an absolutely continuous map, i.e.\ the \( L^{\infty}
\) norm of the solution is bounded by the \( L^{1} \) norm of the
integrand.

The resulting analytic map
\[ \mathbb{S}^{1}\to C^{0}(\mathbb{S}^{2},\mathrm{U}(N)), \]
can be continued to an analytic map
\begin{displaymath}
 \{1-\epsilon<\abs{\lambda}<1+\epsilon\}\to
 C^{0}(\mathbb{S}^{2},\mathrm{GL}(N)),
\end{displaymath}
on some annulus containing the equator.  Since \eqref{Du-Phi} is the
`real form' of the `holomorphic' equation \eqref{op}, this solution
defines a global trivialisation of the bundle \( \mathcal{V} \) on a
deleted neighbourhood of the equator, and we can use it to define the
holomorphic structure of the bundle over the equator.  Grauert's
Theorem implies
that the bundle is trivial on generic fibres.

To see  this
rigorously, observe that  \eqref{op} and \eqref{Du-Phi} can both be
completed to the system \eqref{zccond} by adding a second operator
which has nonzero \( \frac{\partial}{\partial t} \) component.  The
solution to  \eqref{Du-Phi} has a unique extention to a neighbourhood
of \( \{t=t_{0}\} \) and the extension is
 in the kernel of this second operator.  The
resulting solution is a solution to
 \eqref{zccond} and hence a solution to \eqref{op}.
The important point is that null \eqref{A-3} monodromy insures that the
solution is defined on the compactification of  \( \{t=t_{0}\} \)
to a sphere, otherwise the resulting holomorphic trivialisation would
have been for a neighbourhood in \( T\mathbb{P}^{1} \) and not in
\( \T \).

\setcounter{subsection}{\value{equation}}
\addtocounter{equation}{1}
\subsection{Reality}
Reality of the associated bundle is independent of the boundary
conditions and gauge fixing, and is implied by the analogous property
for arbitrary solutions of the Bogomolny equations.
The simplest way to see it in this case is via the formula
$$\lambda^{2}\overline{\tilde{f}^{-1}\sigma^{*}\overline{\nabla}_{\eta}
\tilde{f}\tilde{f}^{-1}}^{t}=\overline{\nabla}_{\eta}\overline{(\sigma^{*}
\tilde{f})}^{t}$$
for a local gauge, \( \tilde{f} \), which shows that holomorphic
gauges are transformed into antiholomorphic gauges of the dual bundle.

\setcounter{subsection}{\value{equation}}
\addtocounter{equation}{1}
\subsection{The section at infinity and the framing}
Over a (possibly pinched) tubular neighbourhood of \( G_{\infty} \),
the section at infinity, the iterative Cauchy-kernel argument
 defines a holomorphic framing.  The radius of the tubular
neighbourhood depends on an energy estimate and is nonzero away from
the equator.
Since the data are holomorphic in \( \lambda \), the result is
holomorphic in base and fibre directions and on \( G_{\infty} \)
agrees with the trivialisation coming from integrating \eqref{A-3} from
infinity.  The resulting trivialization of
\( \mathcal{V}|_{G_{\infty}} \) defines the canonical framing.
Grothendieck's theorem on
formal functions implies that any bundle trivial on a
rational curve of negative self-intersection is trivial on a
neighbourhood of the curve.  So the bundle is actually trivial on a
neighbourhood of \( G_{\infty} \).

\section{Inverse construction : compact twistor fibration}

The inverse construction follows the inverse construction
of \( \nabla,\Phi \) due to Hitchin.  To accommodate the boundary
condition, we need to extend the twistor fibration (and definition of \( J
\)) to a compact twistor fibration.

The first step is to embed $T\mathbb{P}^{1}$ as the nonsingular part of
the singular quadric
$Q\deq\{\beta^{2}=\alpha\gamma\}\subset\mathbb{P}^{3}$ by
$$(\lambda,\eta)\mapsto [1,-2i\lambda,-\lambda^{2},-\eta]=
[\alpha,\beta,\gamma,\delta]$$
(in terms of affine coordinates $\eta\frac{d}{d\lambda}\in
T\mathbb{P}^{1}$ and homogeneous coordinates on $\mathbb{P}^{3}$).
Since the bundle is trivial on a (complex) neighbourhood of the
section at infinity, $\mathcal{V}$ pushes
down via the collapsing map $\T\to Q$
($G_{\infty} \to$ singular point) to a bundle on $Q$.

The next step is to construct the compact double twistor fibration:
\begin{gather*}
	\begin{aligned}
	X\deq&\left\{
	\begin{matrix}
		a\alpha+b\beta+  c\gamma+d\delta=0 \\
		\beta^{2}=\alpha\gamma
	\end{matrix}
	 \right\}
	  \subset\mathbb{P}^{3} \times \mathbb{P}^{3}\\
    &\pi_{1}\swarrow \qquad\quad  \qquad\quad \searrow \pi_{2}\\
	\mathbb{R}^{2+1}\subset&\mathbb{P}^{3}\qquad \qquad \qquad \qquad
	\ \ Q
	\end{aligned}
\end{gather*}
Grauert's Theorem implies that pulling $\mathcal{V}$ back
to $X$
and pushing it forward to $\mathbb{P}^{3}$ gives a coherent sheaf
which we assume
is locally-free on a neighbourhood of $\mathbb{R}^{2+1}\subset
\mathbb{C}^{3}\subset \mathbb{P}^{3}$.
(We will show in a future paper
 that this assumption is unneccessary, i.e.\ that real
bundles which are trivial on equatorial fibres are necessarily
trivial on real sections.)
Call the new sheaf \(
\mathcal{W}\to\mathbb{P}^{3} \).
Fixing a fibre
$P_{\lambda}\subset\T$ such that
$\mathcal{V}|_{P_{\lambda}}$ is trivial,  the composition
$$\mathcal{W}_{y}=H^{0}(G_{y},\mathcal{V})\overset{\text{eval}}{\cong}
\mathcal{V}|_{G_{y}\cap
P_{\lambda}}\overset{\text{eval}}{\cong}H^{0}(P_{\lambda},
\mathcal{V})\cong\mathbb{C}^{N},$$
where \( G_{y}\deq \pi_{2*}\pi_{1}^{-1}(y) \),
gives a natural frame of $\mathcal{W}|_{Y}$,
\[ Y\deq\{y\in\mathbb{P}^{3}:(\pi_{2}^{*}\mathcal{V})|_{\pi_{1}^{-1}(y)}
\text{ is trivial}\}.
 \]
In particular, the standard gauge comes from the fixed framing of
$\mathcal{V}|_{P_{-1}}$, and \( J \) is the
gauge transformation from the \( P_{-1} \) to the \( P_{1} \) framing.
It follows that \( J \)  extends
meromorphically to \( \mathbb{P}^{3} \).

In terms of projective coordinates $[a,b,c,d]$ on
$\mathbb{P}^{3}$, the `finite' hyperplane sections
$\{[a,b,c,1]\}=\mathbb{C}^{3}\subset\mathbb{P}^{3}$ represent the sections
$\{\eta=a-2ib\lambda-c\lambda^{2}\}$ of $T\mathbb{P}^{1}$.  The
`infinite'
hyperplanes $\{[a,b,c,0]\}$ represent the completion of the linear
system on $\T$ to include the family of
divisors $G_{[a,b,c,0]}\deq
G_{\infty}+P_{\lambda_{0}}+P_{\lambda_{-{1}}}$ (where
$a-2ib\lambda_{i}-c\lambda_{i}^{2}=0$).  We know that the set of such
hyperplane sections over which $\mathcal{V}$ is trivial is open and
includes the circle
$\{G_{\infty}+2P_{\lambda}:\lambda\in\mathbb{S}^{1}\}$.  The
intersection $G_{[a,b,c,0]}\cap P_{\lambda}$ is either $P_{\lambda}\cap
G_{\infty}\text{ or } P_{\lambda}$.  Since $P_{\lambda}$ was taken so
that $\mathcal{V}|_{P_{\lambda}}$ is trivial, the definition of the
standard and $P_{\lambda}$ frames extends to an open set of points of
the plane at infinity in $\mathbb{P}^{3}$, and they agree
on this set by definition.  In particular, $J_{\theta}$, the transformation
from the $P_{-1}$ frame to the $P_{e^{i\theta}}$ frame is the identity on the
infinite points.  Since \( J_{e^{i\theta}} \) is in the kernel of
\eqref{A-3} and is defined on compactified space planes, \eqref{A-3} has null
monodromy.

Since \( J \) is analytic by construction, we can use
power series:
Let $b/a$, $c/a$, $d/a$ be affine coordinates on $\mathbb{P}^{3}$
centred at a point at infinity.  $J$ is defined on an open set in
this coordinate chart containing $(0,1,0)$.  The plane at infinity is
cut out by the equation $d/a=0$.  Since $J|_{\{d/a=0\}}=\mathbb{I}$,
we can expand $J$ in a
power series
\begin{equation}\label{powerseries}
\begin{aligned}
	J=\mathbb{I}&+\sum_{\substack{i\geq 1\\ j,k\geq 0}}
	\left(\frac{d}{a}\right)^{i}\left(\frac{b}{a}\right)^{j}\left(\frac{c}{a}
	\right)^{k}J_{ijk}\\
	=\mathbb{I}&+\left(\frac{d}{a}\right)\sum_{k \geq
	0}^{}\left(\frac{c}{a}\right)^{k}J_{10k}\\
	&+\left(\frac{d}{a}
	\right)^{2} \left(\frac{b}{d}\right)\sum_{\substack{j \geq 1\\k \geq 0}}^{}
	\left(\frac{b}{a}\right)^{j-1}\left(\frac{c}{a}\right)^{k}
	J_{ijk}\\
	&+ \left(\frac{d}{a}\right)^{2}\sum_{\substack{i \geq 2\\j,k \geq0}}
	\left(\frac{d}{a}\right)^{i-2}\left(\frac{b}{a}\right)^{j}
	\left(\frac{c}{a}\right)^{k}J_{ijk}\\
		=\mathbb{I}&
	+1/rJ_{1}(\theta)+1/r^{2}J_{2}(\theta,t)+1/r^{2}J_{3}(\theta,1/r,t)
\end{aligned}
\end{equation}
where we have used $d/a=1/z=1/r e^{-i\theta}$, $b/a=-2i t/r
e^{i\theta}$, $c/a=e^{-2i\theta}$ in terms of cylindrical coordinates
on $\mathbb{C}^{3}$, which shows that $J$ satisfies the required
boundary conditions \eqref{A-2}.

This completes the proof that solutions of Ward's equations satisfying
the boundary conditions \eqref{A-2} with null \eqref{A-3} monodromy
are in one to one correspondence with
framed holomorphic bundles over $\T$ which
satisfy a reality and certain triviality conditions.

\begin{remark} \rm
	In a future paper, we will use monads to show that triviality on
equatorial
	fibres plus reality implies triviality on real fibres.
\end{remark}

\begin{remark} \rm
	It follows from \eqref{powerseries} that the energy decays as \(
	\frac{1}{r^{4}} \) as \( r\to\infty \),
	as Ward observed for his solutions.
	This is a property of analytic functions on \( \mathbb{S}^{2}\times
\mathbb{R} \) which
	are constant on \( \{\infty\}\times\mathbb{R} \).
\end{remark}

\newpage
\bibliographystyle{amsalpha}

\end{document}